\documentclass[oneside,reqno]{amsart}

\usepackage{latexsym}
\usepackage{amssymb}
\usepackage{amsmath,amsfonts,amssymb,euscript,array,enumerate,amsfonts,mathrsfs}
\usepackage{graphicx}
\usepackage{footnote}

\newtheorem{Theorem}{Theorem}[section]
\newtheorem{Definition}{Definition}[section]
\newtheorem{Proposition}{Proposition}[section]

\newtheorem{Remark}{Remark}[section]

\def\esssup_#1{\underset{#1}{\mathrm{ess\,sup\, }}}
\def\essinf_#1{\underset{#1}{\mathrm{ess\,inf\, }}}

\def \N{\mathbb{N}}
\def \R{\mathbb{R}}

\def \E{\mathbb{E}}
\def \F{\mathbb{F}}

\def \P{\mathbb{P}}
\def \Q{\mathbb{Q}}

\def \S{\mathbb{S}}

\def \ep{\hbox{ }\hfill$\Box$}

\def\reff#1{{\rm(\ref{#1})}}

\def\beqs{\begin{eqnarray*}}
\def\enqs{\end{eqnarray*}}
\def\beq{\begin{eqnarray}}
\def\enq{\end{eqnarray}}

\begin{document}

\title[Viscosity Solutions of PPDEs and FBSDEs]{Viscosity Solutions of Path-Dependent PDEs and Non-Markovian Forward-Backward Stochastic Equations}

\author{Andrea Cosso$^*$}
\thanks{$^*$Dipartimento di Matematica, Politecnico di Milano, piazza Leonardo da Vinci 32, 20133 Milano, Italy}
\curraddr{}
\email{andrea.cosso@polimi.it}

\subjclass[2010]{35D40, 35K10, 60H10, 60H30.}

\keywords{Non-Markovian FBSDEs, path-dependent PDEs, functional It\^o calculus, viscosity solutions.}

\begin{abstract}
It is known that Markovian forward-backward stochastic differential equations provide nonlinear Feynman-Kac representation formulae for semilinear parabolic PDEs. We show that non-Markovian forward-backward stochastic differential equations provide nonlinear Feynman-Kac formulae for semilinear path-dependent PDEs. This extends the result proved in Ekren, Keller, Touzi, and Zhang \cite{ektz} to the case with a possibly degenerate diffusion coefficient in the forward dynamics.
\end{abstract}

\maketitle

\section{Introduction}

It is well-known, since Peng~\cite{peng91} and Pardoux and Peng~\cite{parpen92}, that Markovian forward-backward stochastic differential equations (FBSDEs) provide a probabilistic representation (also called nonlinear Feynman-Kac formula) for semilinear parabolic partial differential equations.

Recently, Dupire~\cite{dupire} gave rise to a new branch of stochastic calculus, namely functional It\^o calculus. This latter, which has been rigorously developed by Cont and Fourni\'{e}~\cite{cf13}, led also to a generalized family of partial differential equations, known as path-dependent PDEs (PPDEs).

The problem of providing a definition of weak solution (in particular, viscosity solution) for path-dependent PDEs has attracted a great interest. Especially, because of the hypotheses which guarantee the existence of classical solutions are quite strong, see Peng and Wang~\cite{pengwang}. In addition, it seems reasonable that non-Markovian forward-backward stochastic differential equations should provide a ``weak'' solution to semilinear path-dependent PDEs, whenever a well-posedness result for non-Markovian FBSDEs is in force.

Ekren et al.~\cite{ektz} recently provided a definition of viscosity solutions to PPDEs, replacing the classical minimum/maximum property, which appears in the definition of standard viscosity solution, with an optimal stopping problem under nonlinear expectation. In this paper we adopt the definition of viscosity solution designed by Ekren, Touzi, and Zhang in \cite{etzI}, which refines and generalizes to fully nonlinear parabolic PPDEs the definition for the semilinear case presented in \cite{ektz}. However, we note that other definitions of viscosity solutions for PPDEs were given by Peng~\cite{peng12} and Tang and Zhang~\cite{tangzhang}.

In the papers \cite{ektz,etzI} and also in the paper by Henry-Labord\`ere, Tan, and Touzi~\cite{h-ltt13}, it is shown that non-Markovian FBSDEs provide a viscosity solution to semilinear parabolic path-dependent PDEs, with the additional requirements that the driving Brownian motion and the forward process take values in the same Euclidean space and the diffusion part is nondegenerate. In other words, the diffusion coefficient $\sigma$ in the forward dynamics is a square matrix such that $\sigma\sigma^T>0$. Our aim is to show that their result holds even in the degenerate case with $\sigma$ not necessarily square (this is in particular emphasized in Remark~\ref{R:Enlarge}). We also exploit the comparison Theorem 3.4 in \cite{etzII} to give a well-posedness result for our path-dependent PDE. However, even if Theorem 3.4 in \cite{etzII} is, up to now, the finest comparison theorem for parabolic path-dependent PDEs, it requires strong assumptions on the forward coefficients, e.g., the nondegeneracy of the diffusion coefficient. We also mention that another possibility would be to adapt Theorem 7.4 in \cite{etzI} to our context. Indeed, Theorem 7.4 may be seen as a sharper result which can be obtained when we have a representation formula for a viscosity solution to our path-dependent PDE. Nevertheless, Theorem 7.4 uses in a crucial way regularity theory for classical semilinear parabolic PDEs. Hence we need again the nondegeneracy of the diffusion part, unless we use different techniques to regularize viscosity solutions. We leave this issue for future research.

We give an outline of the problem and we also provide a way to deduce formally the form of the PPDE starting from the non-Markovian FBSDE. On the Wiener space, for every $0 \leq t \leq T$ and $x\in\mathbb{R}^{n}$, consider the following non-Markovian forward-backward stochastic differential equation:
\[
\hspace{-0.05cm}\begin{cases}
X_s=x+\int_{t}^{s}b(r,X_r)\textup{d}r+\int_{t}^{s}\sigma(r,X_r)\textup{d}B_r \\
Y_s=g(X_T)+\int_{s}^{T}f(r,X_r,Y_r,Z_r)\textup{d}r-\int_{s}^{T}Z_r\textup{d}B_r
\end{cases}
\]
where $b$, $\sigma$, $f$, and $g$ are stochastic. Suppose that the above FBSDE admits a unique solution $(X^{t,x},Y^{t,x},Z^{t,x})=\{(X_s^{t,x},Y_s^{t,x},Z_s^{t,x}),\,t \leq s \leq T\}$. Let us introduce the function
\[
v(t,x):=Y_t^{t,x}, \qquad \forall(t,x)\in[0,T]\times\mathbb{R}^{n}.
\]
Observe that the function $v$ depends also on $\omega$, therefore it is more properly a random field. This random field is formally related to the following non-Markovian backward stochastic partial differential equation (BSPDE):
\[
\begin{cases}
\!\text{d}v(t,x)=-\Big(\!\mathcal{L}v(t,x) + f\big(t,x,v(t,x),u(t,x)+D_{x}v(t,x)\sigma(t,x)\big) + \text{tr}\big(D_{x}u(t,x)\sigma(t,x)\big)\!\Big)\text{d}t \\
\qquad\qquad\;\;+\,u(t,x)\text{d}B_t, \hspace{6cm} (t,x)\in[0,T)\times\mathbb{R}^{n}\\
\!v(T,x)=g(x), \hspace{7.9cm}x\in\mathbb{R}^{n}
\end{cases}
\]
where the operator $\mathcal{L}$ is given by
\[
\mathcal{L}v=\langle b, D_{x}v\rangle+\frac{1}{2}\text{tr}\big(D_{xx}^2 v \, \sigma\sigma'\big)
\]
and
\[
u(t,x):=Z_t^{t,x}-D_{x}v(t,x)\sigma(t,x), \qquad \forall(t,x)\in[0,T]\times\mathbb{R}^{n}.
\]
We refer to Ma, Yin, and Zhang \cite{myz12} for recent well-posedness results for BSPDEs.

Our aim, instead, is to prove that the random field $v$ is a viscosity solution to a path-dependent PDE. To derive the expression of this PPDE, a possibility consists in using the results presented in \cite{ektz,etzI,pengwang}, where starting from a non-Markovian BSDE the authors deduce the form of the corresponding path-dependent PDE. We can exploit these results in our context, only formally, to associate a path-dependent PDE to our BSPDE. Indeed, for every fixed $x\in\mathbb{R}^{n}$, we can think of our non-Markovian BSPDE as a non-Markovian BSDE, then we deduce the PPDE for $v$:
\[
\begin{cases}
\partial_{t}v(t,x)+\mathcal{L}v(t,x)+\text{tr}\big(D_{x}\partial_\omega v(t,x)\sigma(t,x)\big)+\frac{1}{2}\text{tr}\big(\partial_{\omega\omega}^2v(t,x)\big)\\
\hspace{3.22cm}+\,f\big(t,x,v(t,x),\partial_\omega v(t,x)+D_{x}v(t,x)\sigma(t,x)\big)=0, \;\; (t,x)\in[0,T)\times\mathbb{R}^{n} \\
v(T,x)=g(x), \hspace{9.2cm}x\in\mathbb{R}^{n}
\end{cases}
\]
where $\partial_{t}v$, $\partial_\omega v$, and $\partial_{\omega\omega}^2 v$ are the so-called pathwise or functional derivatives. We note that the above equation is not precisely a path-dependent PDE, due to the presence of the classical derivatives $D_{x}v$ and $D_{xx}^{2}v$. For this reason, in the paper, we refer to it as mixed-path-dependent PDE (MPPDE). Finally, it is simple to associate a path-dependent PDE to the above mixed-path-dependent PDE, simply ``replacing'' all classical derivatives with pathwise derivatives. Our task is now to prove rigorously that $v$ is related to this path-dependent PDE.

The plan of the paper is as follows: section \ref{S:Preliminaries} is devoted to fix the notations and to present some fundamental tools of functional It\^o calculus. In section \ref{S:FBSDE} we introduce the non-Markovian forward-backward equation and we prove the boundedness and the uniform continuity of the value function associated to this FBSDE. Finally, in section \ref{S:PPDE} we study the path-dependent PDE related to the non-Markovian FBSDE. More precisely, we prove, using the method of enlarging the space, that the value function is a viscosity solution to the PPDE and we state a well-posedness result.

\section{Preliminaries}
\label{S:Preliminaries}

In this section we introduce the main tools of functional It\^o calculus needed later, following closely \cite{etzI,etzII}.

\subsection{Non-anticipative functionals on continuous paths}

For any $0 \leq t \leq T$, $k\in\N\backslash\{0\}$, let $\Omega^{t,k}:=\{\omega\in C([t,T],\R^k)\colon \omega_t = 0\}$ be the canonical space of $\R^k$-valued continuous paths on $[t,T]$ starting from the origin, equipped with the uniform norm $\|\omega\|_T^{t,k} := \sup_{t \leq s \leq T} |\omega_s|$, $\omega \in \Omega^{t,k}$. Let $B^{t,k}$ denote the canonical process on $\Omega^{t,k}$, $\F^{t,k}=(\mathcal{F}_s^{t,k})_{t \leq s \leq T}$ denote the natural filtration generated by $B^{t,k}$, and $\P_0^{t,k}$ the Wiener measure on $(\Omega^{t,k},\mathcal{F}_T^{t,k})$. Now, define $\Lambda^{t,k} := [t,T]\times\Omega^{t,k}$ and $\textbf{d}_\infty^{t,k}\colon\Lambda^{t,k}\times\Lambda^{t,k}\rightarrow[0,+\infty)$ as follows:
\[
\textbf{d}_\infty^{t,k} \big((s,\omega),(s',\omega')\big) := |s-s'| + \|\omega_{\cdot \wedge s} - \omega'_{\cdot \wedge s'}\|_T^{t,k},
\]
for all $(s,\omega),(s',\omega')\in\Lambda^{t,k}$. For $\Omega^{t,k}$, $B^{t,k}$, $\F^{t,k}$, $\mathcal{F}_s^{t,k}$, $\P_0^{t,k}$, $\|\cdot\|_T^{t,k}$, $\Lambda^{t,k}$, and $\textbf{d}_\infty^{t,k}$, we drop the superscript $t$ whenever $t=0$.

\vspace{3mm}

We introduce the following spaces of non-anticipative functionals on continuous paths. We denote by $\S^k$ the set of $k\times k$ symmetric matrices.

\begin{Definition}
\label{D:Consistent}
Let $u:\Lambda^{t,k}\rightarrow\R$ be $\F^{t,k}$-progressively measurable.
\begin{itemize}
\item[(i)] We say $u\in C^0(\Lambda^{t,k})$ $($resp., $u\in UC(\Lambda^{t,k})$$)$ if $u$ is continuous $($resp., uniformly continuous$)$ in $(s,\omega)$ under ${\bf d}_\infty^{t,k}$.
\item[(ii)] We say $u\in C_b^0(\Lambda^{t,k})$ $($resp., $u\in UC_b(\Lambda^{t,k})$$)$ if $u\in C^0(\Lambda^{t,k})$ $($resp., $u\in UC(\Lambda^{t,k})$$)$ and $u$ is bounded.
\item[(iii)] We say $u=(u_i)_{i=1,\ldots,k}\in C^0(\Lambda^{t,k},\R^k)$ $($resp., $u=(u_{i,j})_{i,j=1,\ldots,k}\in C^0(\Lambda^{t,k},\S^k)$$)$ if each component $u_i\in C^0(\Lambda^{t,k})$ $($resp., $u_{i,j}\in C^0(\Lambda^{t,k})$$)$. In an obvious way we introduce similar notations for $C_b^0$, $UC$, and $UC_b$.
\end{itemize}
\end{Definition}

\begin{Definition}
We denote by $\underline{\mathcal{U}}^k$ the collection of all $\F^k$-progressively measurable stochastic processes $u\colon\Lambda^k\rightarrow\R$ such that:
\begin{itemize}
\item[(i)] $u$ is bounded from above.
\item[(ii)] For all $\omega\in\Omega^k$, $t \mapsto u(t,\omega)$ is c\`adl\`ag with positive jumps.
\item[(iii)] There exists a modulus of continuity $\rho$ such that
\[
u(t,\omega) - u(t',\omega') \leq \rho\big(\textup{\textbf{d}}_\infty^k((t,\omega),(t',\omega'))\big),
\]
for all $(t,\omega),(t',\omega')\in\Lambda^k$, with $t \leq t'$.
\end{itemize}
We denote by $\overline{\mathcal{U}}^k$ the collection of all stochastic processes $u$ such that $-u \in \underline{\mathcal{U}}^k$.
\end{Definition}

Now we introduce the concatenation operator. Let $0\leq t\leq s \leq T$. For any $\omega\in\Omega^{t,k}$ and $\tilde\omega\in\Omega^{s,k}$, we define the concatenation of $\omega$ and $\tilde\omega$ at $s$ as:
\[
(\omega\otimes_s\tilde\omega)(r) := \omega_r 1_{[t,s)}(r) + (\omega_{s} + \tilde\omega_r) 1_{[s,T]}(r), \qquad t \leq r \leq T.
\]
Note that $\omega\otimes_s\tilde\omega$ lies in $\Omega^{t,k}$. Then, for any $\mathcal{F}_T^{t,k}$-measurable random variable $\xi\colon\Omega^{t,k}\rightarrow\R$ and $\F^{t,k}$-progressively measurable stochastic process $u\colon[t,T]\times\Omega^{t,k}\rightarrow\R$, we define the shifted $\mathcal{F}_T^{s,k}$-measurable random variable $\xi^{s,\omega}\colon\Omega^{s,k}\rightarrow\R$ and $\F^{s,k}$-progressively measurable stochastic process $u^{s,\omega}\colon[s,T]\times\Omega^{s,k}\rightarrow\R$ as:
\[
\xi^{s,\omega}(\tilde\omega) := \xi(\omega\otimes_s\tilde\omega), \qquad u_r^{s,\omega}(\tilde\omega) := u_r(\omega\otimes_s\tilde\omega),
\]
for all $\tilde\omega\in\Omega^{s,k}$, $s \leq r \leq T$.

Finally, let $0 \leq t \leq T$ and ${\mathcal{T}}^{t,k}$ denote the set of all $\F^{t,k}$-stopping times on $[t,T]$. Then we define ${\mathcal{H}}^{t,k}\subset{\mathcal{T}}^{t,k}$ as the subset of hitting times $\text{\footnotesize H}$ of the form:
\[
\text{\footnotesize H} := \inf\big\{s\in[t,T]\colon B_s^{t,k} \in O^c\big\}\wedge t_0,
\]
for some $t < t_0 \leq T$ and some open and convex set $O \subset \R^k$ containing the origin. For ${\mathcal{T}}^{t,k}$ and ${\mathcal{H}}^{t,k}$ we drop the superscript $t$ whenever $t=0$.

\subsection{Nonlinear expectation and pathwise derivatives}

For all $L>0$ and $0 \leq t \leq T$, let ${\mathcal{P}}_L^{t,k}$ denote the set of probability measures $\P^{t,k}$ on $(\Omega^{t,k},\mathcal{F}_T^{t,k})$ such that there exist an $\F^{t,k}$-progressively measurable $\R^k\times\S^k$-valued process $(\alpha,\beta)$ and a $k$-dimensional $\P^{t,k}$-Wiener process $\tilde B$ satisfying:
\[
|\alpha_s| \leq L, \qquad |\beta_s| \leq L, \qquad dB_s^{t,k} = \alpha_s ds + \beta_s d\tilde B_s, \quad t \leq s \leq T, \; \P^{t,k}\text{-a.s.}
\]
Being $\beta$ an $\R^{k\times k}$-valued stochastic process, when we write $|\beta_s|$ we suppose that $|\cdot|$ is a norm on $\R^{k\times k}$ fixed throughout the paper. Clearly, the same remark applies to the norm $|\cdot|$ on $\R^k$ applied to $\alpha_s$.

\begin{Remark}
\label{BetaSymm}
{\rm Note that it is not a restriction to consider $\beta$ symmetric. Indeed, if $\beta$ is an $\F^{t,k}$-progressively measurable $\R^{k \times k}$-valued process such that $|\beta_s| \leq L$ for all $t \leq s \leq T$, $\P^{t,k}$-a.s., then $\int_t^s\beta_r d\tilde B_r$ and $\int_t^s(\beta_r\beta_r^*)^{1/2}d\tilde B_r$ have the same distribution.
}
\end{Remark}

Now, let $\xi\colon\Omega^{t,k}\rightarrow\R$ (resp., $\eta\colon\Omega^{t,k}\rightarrow\R$) be a bounded from above (resp., from below) $\mathcal{F}_T^{t,k}$-measurable random variable. Then, for all $0 \leq t \leq T$, we introduce the \emph{nonlinear expectations}:
\begin{equation}
\label{nonlinear_expectations}
\underline{\mathcal{E}}_t^L [\xi] := \inf_{\P^{t,k}\in{\mathcal{P}}_L^{t,k}}\E^{\P^{t,k}}[\xi], \qquad\qquad
\overline{\mathcal{E}}_t^L [\eta] := \sup_{\P^{t,k}\in{\mathcal{P}}_L^{t,k}}\E^{\P^{t,k}}[\eta].
\end{equation}
Let $Y\colon[0,T]\times\Omega^k\rightarrow\R$ be a bounded $\F^k$-progressively measurable stochastic process. Then, for all $(t,\omega)\in\Lambda^k$, we introduce the \emph{nonlinear Snell envelopes}:
\begin{equation}
\label{nonlinear_stopping}
\underline{\mathcal{S}}_t^L [Y](\omega) := \inf_{\tau\in{\mathcal{T}}^{t,k}}\underline{\mathcal{E}}_t^L[Y_\tau^{t,\omega}], \qquad\qquad
\overline{\mathcal{S}}_t^L [Y](\omega) := \sup_{\tau\in{\mathcal{T}}^{t,k}}\overline{\mathcal{E}}_t^L[Y_\tau^{t,\omega}].
\end{equation}
These nonlinear optimal stopping problems, deeply studied in \cite{etzOptStop}, are essential for the definition of viscosity solution of path-dependent PDEs.

We conclude this section by introducing the pathwise derivatives. Let $u\colon\Lambda^{t,k}\rightarrow\R$ be $\F^{t,k}$-progressively measurable. Firstly we define the \emph{pathwise time derivative} as follows:
\[
\partial_t u(s,\omega):=\lim_{h\rightarrow0^+} \frac{ u(s+h,\omega_{\cdot \wedge s}) - u(s,\omega)}{h}, \qquad s<T
\]
and
\[
\partial_t u(T,\omega):=\lim_{s\rightarrow T^-}\partial_t u(s,\omega), \qquad s=T.
\]
We define the \emph{pathwise spatial derivatives} via the functional It\^o's formula, as in \cite{etzI,etzII}. To this end, we denote:
\[
{\mathcal{P}}_\infty^{t,k} := \bigcup_{L>0} {\mathcal{P}}_L^{t,k}, \qquad\qquad 0 \leq t \leq T.
\]

\begin{Definition}
We say $u\in C_b^{1,2}(\Lambda^{t,k})$ if $u\in C_b^0(\Lambda^{t,k})$, $\partial_t u\in C_b^0(\Lambda^{t,k})$, and there exist $\partial_\omega u \in C_b^0(\Lambda^{t,k},\R^k)$, $\partial_{\omega\omega}^2 u \in C_b^0(\Lambda^{t,k},\S^k)$ such that for any $(s,\omega) \in [t,T)\times\Omega^{t,k}$ and any $\P^{s,k} \in {\mathcal{P}}_\infty^{s,k}$, $u^{s,\omega}$ is a local $\P^{s,k}$-semimartingale and the \textbf{functional It\^o's formula} holds true, $\P^{s,k}$-a.s.,
\[
du^{s,\omega}(r,B^{s,k}) = \partial_t u^{s,\omega}(r,B^{s,k}) dr + \frac{1}{2}\textup{tr}\big(\partial_{\omega\omega}^2 u^{s,\omega}(r,B^{s,k}) d\langle B^{s,k}\rangle_r\big) + \partial_\omega u^{s,\omega}(r,B^{s,k}) dB_r^{s,k}
\]
for all $s \leq r \leq T$.
\end{Definition}

For ${\mathcal{P}}_L^{t,k}$, ${\mathcal{P}}_\infty^{t,k}$, $\P^{t,k}$, $\underline{\mathcal{E}}_t^L$, $\overline{\mathcal{E}}_t^L$, $\underline{\mathcal{S}}_t^L$, and $\overline{\mathcal{S}}_t^L$, we drop the $t$ whenever $t=0$.

\subsection{Regular conditional probability distribution}
\label{rcpd}

We recall here some standard facts about regular conditional probability distributions for reader's convenience. For more details, see Section 1.3 in \cite{sv06}, Section 4 in \cite{stz13}, and the Appendix in \cite{bh}.

Let $\P$ be an arbitrary probability measure on $(\Omega^k,\mathcal{F}_T^k)$ and $\tau$ be an $\F^k$-stopping time on $[0,T]$. Then, in light of Theorem 1.3.4 in \cite{sv06}, there exists a regular conditional probability distribution (r.c.p.d.) $(\Q_\tau^\omega)_{\omega\in\Omega^k}$ of $\P$ given $\mathcal{F}_\tau^k$, namely:
\begin{itemize}
\item[(i)] For each $\omega\in\Omega^k$, $\Q_\tau^\omega$ is a probability measure on $(\Omega^k,\mathcal{F}_T^k)$;
\item[(ii)] For each $A\in\mathcal{F}_T^k$, the mapping $\omega\mapsto\Q_\tau^\omega(A)$ is $\mathcal{F}_\tau^k$-measurable;
\item[(iii)] For each bounded $\mathcal{F}_T^k$-measurable random variable $\xi$, it holds for $\P$-a.e. $\omega\in\Omega^k$ that $\E^\P(\xi|\mathcal{F}_\tau^k)(\omega)=\E^{\Q_\tau^\omega}(\xi)$;
\item[(iv)] For each $\omega\in\Omega^k$, $\Q_\tau^\omega(\omega\otimes_{\tau(\omega)}\Omega^{\tau(\omega),k})=1$.
\end{itemize}
By property (iv) above, for any fixed $\omega\in\Omega^k$ we can define a probability measure $\Q^{\tau,\omega}$ on $(\Omega^{\tau(\omega),k},\mathcal{F}_T^{\tau(\omega),k})$ by
\[
\Q^{\tau,\omega}(A) := \Q_\tau^\omega\big(\omega\otimes_{\tau(\omega)}A\big), \qquad \forall\,A\in\mathcal{F}_T^{\tau(\omega),k}.
\]
Then, combining properties (iii) and (iv) above and thanks to a monotone class argument, we have: for any bounded $\mathcal{F}_T^k$-measurable random variable $\xi$, it holds, for $\P$-a.e. $\omega\in\Omega^k$,
\[
\E^\P(\xi|\mathcal{F}_\tau^k)(\omega) = \E^{\Q_\tau^\omega}\big(1_{\{\omega\otimes_{\tau(\omega)}\Omega^{\tau(\omega),k}\}}\xi\big) = \E^{\Q^{\tau,\omega}}\big(\xi^{\tau(\omega),\omega}\big).
\]
Note that the r.c.p.d. $(\Q_\tau^\omega)_{\omega\in\Omega^k}$ is generally not unique. However, it can be proved (see the Appendix in \cite{bh}) that there exists a particular r.c.p.d. $(\Q_\tau^\omega)_{\omega\in\Omega^k}$ of $\P_0^k$ such that $\Q^{\tau,\omega}=\P_0^{\tau(\omega),k}$. Hence, for any bounded $\mathcal{F}_T^k$-measurable random variable $\xi$, it holds, for $\P$-a.e. $\omega\in\Omega^k$,
\[
\E^\P(\xi|\mathcal{F}_\tau^k)(\omega) = \E^{\P_0^{\tau(\omega),k}}\big(\xi^{\tau(\omega),\omega}\big).
\]

\section{Non-Markovian Forward-Backward Equation}
\label{S:FBSDE}

In the present section we fit the non-Markovian forward-backward equation into the framework of functional It\^o calculus.

Let $d,n\in\N\backslash\{0\}$. Then for any $(t,\omega,x)\in\Lambda^d\times\R^n$ we consider the following non-Markovian forward-backward stochastic differential equation:
\begin{equation}
\label{FBSDE}
\begin{cases}
X_s = x + \int_t^s b^{t,\omega}(r,B^{t,d},X_r)dr + \int_t^s \sigma^{t,\omega}(r,B^{t,d},X_r)dB_r^{t,d} \\
Y_s = g^{t,\omega}(B^{t,d},X_T) + \int_s^T f^{t,\omega}(r,B^{t,d},X_r,Y_r,Z_r)dr - \int_s^T Z_rdB_r^{t,d}
\end{cases}
\end{equation}
for all $t \leq s \leq T$, $\P_0^{t,d}$-a.s..

\vspace{3mm}

We impose the following conditions on the forward coefficients.

\vspace{3mm}

\textbf{(HFC)}
\begin{itemize}
\item[] The drift $b\colon\Lambda^d\times\R^n\rightarrow\R^n$ and the diffusion coefficient $\sigma\colon\Lambda^d\times\R^n\rightarrow\R^{n \times d}$ are bounded, $b(t,\cdot,x)$ and $\sigma(t,\cdot,x)$ are $\F^d$-progressively measurable for any $x\in\R^n$, and there exist a constant $C>0$ and a concave modulus of continuity $\rho$ such that:
\[
|b(t,\omega,x) - b(t',\omega',x')| \leq C\big(\rho(|t-t'| + \|\omega_{\cdot\wedge t}-\omega_{\cdot\wedge t'}'\|_T^d) + |x-x'|\big)
\]
and
\[
|\sigma(t,\omega,x) - \sigma(t',\omega',x')| \leq C\big(\rho(|t-t'| + \|\omega_{\cdot\wedge t}-\omega_{\cdot\wedge t'}'\|_T^d) + |x-x'|\big),
\]
for any $t,t'\in[0,T]$, $\omega,\omega'\in\Omega^d$, and $x,x'\in\R^n$.
\end{itemize}

\vspace{3mm}

We make the following assumptions on the BSDE coefficients:

\vspace{3mm}

\textbf{(HBC)}
\begin{itemize}
\item[(i)] \emph{Terminal condition.} $g\colon\Omega^d\times\R^n\rightarrow\R$ is bounded and there exist a constant $C>0$ and a concave modulus of continuity $\rho$ such that:
\[
|g(\omega,x) - g(\omega',x')| \leq C\big(\rho(\|\omega-\omega'\|_T^d) + |x-x'|\big),
\]
for any $\omega,\omega'\in\Omega^d$ and $x,x'\in\R^n$.
\item[(ii)] \emph{Generator.} $f\colon\Lambda^d\times\R^n\times\R\times\R^d\rightarrow\R$ is bounded, $f(t,\cdot,x,y,z)$ is $\F^d$-progressively measurable for any $(x,y,z)\in\R^{n+1+d}$, and there exist a constant $C>0$ and a concave modulus of continuity $\rho$ such that:
\begin{align*}
|f(t,\omega,x,y,z) - f(t',\omega',x',y',z')| \leq C\big(\rho(|t-t'| &+ \|\omega_{\cdot\wedge t}-\omega_{\cdot\wedge t'}'\|_T^d) \\
&+ |x-x'| + |y-y'| + |z-z'|\big),
\end{align*}
for any $t,t'\in[0,T]$, $\omega,\omega'\in\Omega^d$, $x,x'\in\R^n$, $y,y'\in\R$, and $z,z'\in\R^d$.
\end{itemize}

\vspace{3mm}

Under assumptions \textbf{(HFC)} and \textbf{(HBC)} it can be shown, from the standard theory of SDEs and the results presented in \cite{parpen90}, that there exists a unique solution to the above forward-backward stochastic differential equation, denoted by $(X^{t,\omega,x},Y^{t,\omega,x},Z^{t,\omega,x})$ $=$ $\{(X_s^{t,\omega,x},Y_s^{t,\omega,x},Z_s^{t,\omega,x}),\,t \leq s \leq T\}$ such that:
\begin{itemize}
\item[(i)] $X^{t,\omega,x}$ is an $\R^n$-valued $\F^{t,d}$-progressively measurable and continuous process satisfying: for all $p\geq1$ there exists a constant $C_p>0$ such that
\[
\E^{\P_0^{t,d}}\Big[\sup_{t \leq s \leq T}|X_s^{t,\omega,x}|^p\Big] \leq C_p (1 + |x|^p).
\]
\item[(ii)] $Y^{t,\omega,x}$ is a real-valued $\F^{t,d}$-progressively measurable and continuous process satisfying:
\[
\E^{\P_0^{t,d}}\Big[\sup_{t \leq s \leq T}|Y_s^{t,\omega,x}|^2\Big] < \infty.
\]
\item[(iii)] $Z^{t,\omega,x}$ is an $\R^d$-valued $\F^{t,d}$-progressively measurable process satisfying:
\[
\E^{\P_0^{t,d}}\int_t^T|Z_s^{t,\omega,x}|^2ds < \infty.
\]
\end{itemize}

Note that $Y_{t}^{t,\omega,x}$ is $\mathcal{F}_t^{t,d}$-measurable and therefore is a constant. We thus define
\begin{equation}
\label{v}
v(t,\omega,x) := Y_t^{t,\omega,x}, \qquad (t,\omega,x)\in\Lambda^d\times\R^n.
\end{equation}

\begin{Proposition}
\label{UC}
Under assumptions \textup{\textbf{(HFC)}} and \textup{\textbf{(HBC)}}, the map $v$ defined in \reff{v} is bounded and uniformly continuous on $\Lambda^d\times\R^n$.
\end{Proposition}
\emph{Proof.}
We begin by noting that the boundedness of $v$ follows directly from the boundedness of $f$ and $g$. Now we prove the uniform continuity property.

In what follows we shall denote by $C>0$ a generic positive constant depending
only on $T$, the forward coefficients $b,\sigma$, and the backward coefficients $f,g$, which may vary from line to line. Let $(t_1,\omega^1,x_1),(t_2,\omega^2,x_2)\in\Lambda^d\times\R^n$, with $t_1 \leq t_2$. For simplicity of notation, we denote $(X^1,Y^1,Z^1):=(X^{t_1,\omega^1,x_1},Y^{t_1,\omega^1,x_1},Z^{t_1,\omega^1,x_1})$ (resp., $(X^2,Y^2,Z^2):=(X^{t_2,\omega^2,x_2},Y^{t_2,\omega^2,x_2},Z^{t_2,\omega^2,x_2})$) the solution to the FBSDE \reff{FBSDE} starting at time $t=t_1$ (resp., $t=t_2$) from $(\omega,x)=(\omega^1,x_1)$ (resp., $(\omega,x)=(\omega^2,x_2)$).

By taking the conditional expectation $\E_{t_2}^{\P_0^{t_1,d}}$ and using the properties of regular conditional probability distributions (see \cite{sv06} or Section \ref{rcpd}), we deduce that $Y^1$ may be seen as the solution to the following BSDE on $[t_2,T]$: for $\P_0^{t_1,d}$-a.e. $B^{t_1,d}$,
\begin{align*}
Y_t^1 = g^{t_2,\omega^1\otimes_{t_1}B^{t_1,d}}(B^{t_2,d},X_T^1) &- \int_t^T Z_s^1 dB_s^{t_2,d} \\
& + \int_t^T f^{t_2,\omega^1\otimes_{t_1}B^{t_1,d}}(s,B^{t_2,d},X_s^1,Y_s^1,Z_s^1)ds
\end{align*}
for all $t_2 \leq t \leq T$, $\P_0^{t_2,d}$-a.s.. Denote, for $\P_0^{t_1,d}$-a.e. $B^{t_1,d}$,
\[
\bar X_t := X_t^1 - X_t^2, \qquad \bar Y_t := Y_t^1 - Y_t^2, \qquad \bar Z_t := Z_t^1 - Z_t^2
\]
for all $t_2 \leq t \leq T$, $\P_0^{t_2,d}$-a.s.. Then, an application of It\^o's formula to $|\bar Y_t|$ yields, for $\P_0^{t_1,d}$-a.e. $B^{t_1,d}$,
\begin{align*}
|\bar Y_t|^2 = |\bar Y_T|^2 + \int_t^T\bar Y_s\big(\gamma_s + \alpha_s\bar Y_s + \beta_s\bar Z_s\big)ds - \frac{1}{2}\int_t^T|\bar Z_s|^2 ds - \int_t^T \bar Y_s\bar Z_s dB_s^{t_2,d}
\end{align*}
for all $t_2 \leq t \leq T$, $\P_0^{t_2,d}$-a.s.. Here, for $\P_0^{t_1,d}$-a.e. $B^{t_1,d}$, $\alpha$ (resp., $\beta$) is an $\F^{t_2,d}$-progressively measurable $\R$-valued (resp., $\R^d$-valued) stochastic process, such that $|\alpha| \leq C$ (resp., $|\beta| \leq C$). Moreover, for $\P_0^{t_1,d}$-a.e. $B^{t_1,d}$,
\[
\gamma_t := f^{t_2,\omega^1\otimes_{t_1}B^{t_1,d}}(t,B^{t_2,d},X_t^1,Y_t^1,Z_t^1) -f^{t_2,\omega^2}(t,B^{t_2,d},X_t^2,Y_t^1,Z_t^1)
\]
for all $t_2 \leq t \leq T$, $\P_0^{t_2,d}$-a.s.. Now, from Young's inequality and Gronwall's lemma, we get, for $\P_0^{t_1,d}$-a.e. $B^{t_1,d}$,
\begin{equation}
|\bar Y_{t_2}|^2 \leq C\bigg(\E^{\P_0^{t_2,d}}\big[|\bar Y_T|^2\big] + \int_{t_2}^T\E^{\P_0^{t_2,d}}\big[|\gamma_s|^2\big]ds\bigg). \label{estimateY}
\end{equation}
Note that
\begin{align}
|\bar Y_T| &= \big|g^{t_2,\omega^1\otimes_{t_1}B^{t_1,d}}(B^{t_2,d},X_T^1) - g^{t_2,\omega^2}(B^{t_2,d},X_T^2)\big| \notag \\
&\leq C\big(\rho\big(\|\omega^1\otimes_{t_1}B^{t_1,d}\otimes_{t_2}B^{t_2,d} - \omega^2\otimes_{t_2} B^{t_2,d}\|_T^d\big) + |\bar X_T|\big) \notag \\
&\leq C\big(\rho\big(\textbf{d}_\infty^d((t_1,\omega^1),(t_2,\omega^2)) + \|B_{\cdot\wedge t_2}^{t_1,d}\|_T^{t_1,d}\big) + |\bar X_T|\big). \label{estimateg}
\end{align}
Similarly,
\begin{equation}
|\gamma_s| \leq C\big(\rho\big(\textbf{d}_\infty^d((t_1,\omega^1),(t_2,\omega^2)) + \|B_{\cdot\wedge t_2}^{t_1,d}\|_T^{t_1,d}\big) + |\bar X_s|\big). \label{estimatef}
\end{equation}
Now we study the difference $|\bar X_s|=|X_s^1 - X_s^2|$. By taking the conditional expectation $\E_{t_2}^{\P_0^{t_1,d}}$, we deduce that $X^1$ may be seen as the solution to the following forward SDE on $[t_2,T]$: for $\P_0^{t_1,d}$-a.e. $B^{t_1,d}$,
\[
X_t^1 = X_{t_2}^1 + \int_{t_2}^t b^{t_2,\omega^1\otimes_{t_1}B^{t_1,d}}(s,B^{t_2,d},X_s^1)ds + \int_{t_2}^t \sigma^{t_2,\omega^1\otimes_{t_1}B^{t_1,d}}(s,B^{t_2,d},X_s^1)dB_s^{t_2,d},
\]
for all $t_2 \leq t \leq T$, $\P_0^{t_2,d}$-a.s.. Then, using standard arguments, we have: for $\P_0^{t_1,d}$-a.e. $B^{t_1,d}$,
\begin{equation}
\E^{\P_0^{t_2,d}}\big[|\bar X_t|^2\big] \leq C \Big(|\bar X_{t_2}|^2 + \rho\big(\textbf{d}_\infty^d((t_1,\omega^1),(t_2,\omega^2)) + \|B_{\cdot\wedge t_2}^{t_1,d}\|_T^{t_1,d}\big)^2\Big), \label{estimateX}
\end{equation}
for all $t_2 \leq t \leq T$. Hence, from \reff{estimateY}, \reff{estimateg}, \reff{estimatef}, and \reff{estimateX} we obtain: for $\P_0^{t_1,d}$-a.e. $B^{t_1,d}$,
\[
\E^{\P_0^{t_2,d}}\big[|\bar Y_t|^2\big] \leq C \Big(|\bar X_{t_2}|^2 + \rho\big(\textbf{d}_\infty^d((t_1,\omega^1),(t_2,\omega^2)) + \|B_{\cdot\wedge t_2}^{t_1,d}\|_T^{t_1,d}\big)^2\Big),
\]
for all $t_2 \leq t \leq T$. Therefore, for $\P_0^{t_1,d}$-a.e. $B^{t_1,d}$,
\[
|\bar Y_{t_2}| \leq C \Big(|\bar X_{t_2}| + \rho\big(\textbf{d}_\infty^d((t_1,\omega^1),(t_2,\omega^2)) + \|B_{\cdot\wedge t_2}^{t_1,d}\|_T^{t_1,d}\big)\Big).
\]
Thus, noting that $f$ is bounded and $\rho$ is subadditive,
\begin{align}
|v(t_1,\omega^1,x_1) - v(t_2,\omega^2,x_2)| &= |Y_{t_1}^1 - Y_{t_2}^2| \notag \\
&= \Big|\E^{\P_0^{t_1,d}}\Big[Y_{t_2}^1 + \int_{t_1}^{t_2}f^{t_1,\omega^1}(s,B^{t_1},X_s^1,Y_s^1,Z_s^1)ds - Y_{t_2}^2\Big]\Big| \notag \\
&\leq C|t_2 - t_1| + \E^{\P_0^{t_1,d}}\big[|\bar Y_{t_2}|\big] \notag \\
&\leq C|t_2 - t_1| + C\E^{\P_0^{t_1,d}}\Big[|\bar X_{t_2}| + \rho\big(\textbf{d}_\infty^d((t_1,\omega^1),(t_2,\omega^2)) \big) \notag \\
&\quad + \rho\big(\|B_{\cdot\wedge t_2}^{t_1,d}\|_T^{t_1,d}\big)\Big]. \label{estimatev}
\end{align}
Now, note that from an application of It\^o's formula to $X_t^1-x_2$ and Gronwall's lemma, we deduce
\[
\E^{\P_0^{t_1,d}}\big[|\bar X_{t_2}|\big] \leq C\big(\sqrt{|t_1-t_2|} + |x_1-x_2|\big).
\]
Therefore, to prove the uniform continuity of $v$ on $\Lambda^d\times\R^n$, it remains to study the following term: $\E^{\P_0^{t_1,d}}[\rho(\|B_{\cdot\wedge t_2}^{t_1,d}\|_T^{t_1,d})]$. Since $\rho$ is concave, Jensen's inequality yields:
\[
\E^{\P_0^{t_1,d}}\Big[\rho\big(\|B_{\cdot\wedge t_2}^{t_1,d}\|_T^{t_1,d}\big)\Big] \leq \rho\Big(\E^{\P_0^{t_1,d}}\Big[\|B_{\cdot\wedge t_2}^{t_1,d}\|_T^{t_1,d}\Big]\Big) \leq \rho\big(\sqrt{t_2-t_1}\big),
\]
which completes the proof.
\ep

\section{Path-dependent PDE}
\label{S:PPDE}

In the present section we prove that the non-Markovian forward-backward stochastic differential equation is related to a path-dependent PDE. Firstly, it is convenient to \emph{enlarge} the canonical space and to view the FBSDE as defined on $(\Omega^{d+n},\mathcal{F}_T^{d+n},\P_0^{d+n})$. More precisely, for any $0 \leq t \leq T$, let us denote by $\omega^{d+n}=(\omega^d,\omega^n)$, with $\omega^d\in\Omega^{t,d}$ and $\omega^n\in\Omega^{t,n}$, a generic element of $\Omega^{t,d+n}$. Let also $B^{t,d+n}=(B^{t,d},B^{t,n})$ denote the canonical process on $\Omega^{t,d+n}$, being $B^{t,d}$ the first $d$ components and $B^{t,n}$ the last $n$ components of $B^{t,d+n}$. Then for any $(t,\omega^{d+n})\in\Lambda^{d+n}$ we consider the following non-Markovian FBSDE:
\begin{equation}
\label{FBSDEd+n}
\begin{cases}
\hat X_s = \omega_t^n + \int_t^s b^{t,\omega^d}(r,B^{t,d},\hat X_r)dr + \int_t^s \sigma^{t,\omega^d}(r,B^{t,d},\hat X_r)dB_r^{t,d} \\
\hat Y_s = g^{t,\omega^d}(B^{t,d},\hat X_T) + \int_s^T f^{t,\omega^d}(r,B^{t,d},\hat X_r,\hat Y_r,\hat Z_r)dr - \int_s^T \hat Z_rd B_r^{t,d}
\end{cases}
\end{equation}
for all $t \leq s \leq T$, $\P_0^{t,d+n}$-a.s.. Note that under assumptions \textbf{(HFC)} and \textbf{(HBC)} it can be shown that there exists a unique solution to the above forward-backward stochastic differential equation, denoted by $(\hat X^{t,\omega^{d+n}},\hat Y^{t,\omega^{d+n}},\hat Z^{t,\omega^{d+n}})=\{(\hat X_s^{t,\omega^{d+n}},\hat Y_s^{t,\omega^{d+n}},\hat Z_s^{t,\omega^{d+n}}),\,t \leq s \leq T\}$ such that:
\begin{itemize}
\item[(i)] $\hat X^{t,\omega^{d+n}}\colon[t,T]\times\Omega^{t,d+n}\rightarrow\R^n$ is a $\sigma(B_s^{t,d},\,t \leq s \leq T)$-progressively measurable and continuous process satisfying: for all $p\geq1$ there exists a constant $C_p>0$ such that
\[
\E^{\P_0^{t,d+n}}\Big[\sup_{t \leq s \leq T}|\hat X_s^{t,\omega^{d+n}}|^p\Big] \leq C_p (1 + |\omega_t^n|^p).
\]
\item[(ii)] $\hat Y^{t,\omega^{d+n}}\colon[t,T]\times\Omega^{t,d+n}\rightarrow\R$ is a $\sigma(B_s^{t,d},\,t \leq s \leq T)$-progressively measurable and continuous process satisfying:
\[
\E^{\P_0^{t,d+n}}\Big[\sup_{t \leq s \leq T}|\hat Y_s^{t,\omega^{d+n}}|^2\Big] < \infty.
\]
\item[(iii)] $\hat Z^{t,\omega^{d+n}}\colon[t,T]\times\Omega^{t,d+n}\rightarrow\R^d$ is a $\sigma(B_s^{t,d},\,t \leq s \leq T)$-progressively measurable process satisfying:
\[
\E^{\P_0^{t,d+n}}\int_t^T|\hat Z_s^{t,\omega^{d+n}}|^2ds < \infty.
\]
\end{itemize}

Finally, we observe that the two stochastic processes
\[
(\Omega^{t,d},\mathcal{F}_T^{t,d},\F^{t,d},X^{t,\omega^d,\omega_t^n},Y^{t,\omega^d,\omega_t^n},Z^{t,\omega^d,\omega_t^n},\P_0^{t,d})
\]
and
\[
(\Omega^{t,d+n},\mathcal{F}_T^{t,d+n},\F^{t,d+n},\hat X^{t,\omega^{d+n}},\hat Y^{t,\omega^{d+n}},\hat Z^{t,\omega^{d+n}},\P_0^{t,d+n})
\]
have the same law. In particular, the map
\begin{equation}
\label{hatv}
\hat v(t,\omega^{d+n}) := \hat Y_t^{t,\omega^{d+n}}, \qquad (t,\omega^{d+n})\in\Lambda^{d+n}
\end{equation}
is such that
\begin{equation}
\label{vhatv}
\hat v(t,\omega^{d+n}) = v(t,\omega^d,\omega_t^n), \qquad (t,\omega^{d+n})\in\Lambda^{d+n},
\end{equation}
where $v$ is given by \reff{v}.

\vspace{3mm}

To deduce the form of our candidate path-dependent PDE, recall from the introduction that we derived formally a link between the forward-backward stochastic differential equation and a mixed-path-dependent PDE (MPPDE), i.e., a PDE characterized by both classical derivatives and pathwise derivatives. It is then natural to associate to this mixed-path-dependent PDE a path-dependent PDE, as shown below.

Let us consider the following \emph{mixed-path-dependent PDE}:
\begin{align}
\label{MPPDE}
-\partial_t u(t,\omega,x) &- \big\langle b(t,\omega,x),D_x u(t,\omega,x)\big\rangle \\
&\hspace{-1cm} - \frac{1}{2}\text{tr}\Big(D_{xx}^2 u(t,\omega,x)\sigma(t,\omega,x)\sigma^T(t,\omega,x) + 2D_x\partial_\omega u(t,\omega,x)\sigma(t,\omega,x) + \partial_{\omega\omega}^2 u(t,\omega,x)\Big) \notag \\
&\hspace{2.6cm} - f\big(t,\omega,x,u(t,\omega,x),\partial_\omega u(t,\omega,x) + D_x u(t,\omega,x)\sigma(t,\omega,x)\big) = 0, \notag
\end{align}
for all $(t,\omega,x)\in[0,T)\times\Omega^d\times\R^n$, with the \emph{terminal condition}:
\[
u(T,\omega,x) = g(\omega,x), \qquad\qquad (\omega,x)\in\Omega^d\times\R^n.
\]
Here $D_x$ and $D_{xx}$ are the classical derivatives with respect to $x$.

\vspace{3mm}

Now, it is natural to associate to the above mixed-path-dependent PDE the following \emph{path-dependent PDE}:
\begin{align}
\label{PPDE}
-\partial_t \hat u(t,\omega^{d+n}) &- (\hat{\mathcal{L}} \hat u)(t,\omega^{d+n}) \\
&- f\big(t,\omega^d,\omega_t^n,\hat u(t,\omega^{d+n}),\partial_{\omega^d} \hat u(t,\omega^{d+n}) + \partial_{\omega^n} \hat u(t,\omega^{d+n})\sigma(t,\omega^d,\omega_t^n)\big) = 0, \notag
\end{align}
for all $(t,\omega^{d+n})\in[0,T)\times\Omega^{d+n}$, with the \emph{terminal condition}:
\begin{equation}
\label{termcondPPDE}
\hat u(T,\omega^{d+n}) = g(\omega^d,\omega_T^n), \qquad\qquad \omega^{d+n}\in\Omega^{d+n},
\end{equation}
where
\begin{align*}
(\hat {\mathcal{L}} \hat u)(t,\omega^{d+n})&:= \big\langle b(t,\omega^d,\omega_t^n),\partial_{\omega^n} \hat u(t,\omega^{d+n})\big\rangle \\
&\quad + \frac{1}{2}\text{tr}\Big(\partial_{\omega^n\omega^n}^2 \hat u(t,\omega^{d+n}) \sigma(t,\omega^d,\omega_t^n) \sigma^T(t,\omega^d,\omega_t^n) \\
&\quad + 2\partial_{\omega^d\omega^n}^2 \hat u(t,\omega^{d+n}) \sigma(t,\omega^d,\omega_t^n) + \partial_{\omega^d\omega^d}^2 \hat u(t,\omega^{d+n})\Big),
\end{align*}
for all $(t,\omega^{d+n})\in[0,T)\times\Omega^{d+n}$.

Now we give the definition, introduced in \cite{etzI}, of viscosity solution for the path-dependent PDE \reff{PPDE}. Firstly, we define the following sets of test functions.

\begin{Definition}
Let $\hat u \in \underline{\mathcal{U}}^{d+n}$. For every $L>0$ and $(t,\omega^{d+n})\in[0,T)\times\Omega^{d+n}$, define:
\begin{align}
\label{AL}
\underline{\mathcal{A}}^L \hat u(t,\omega^{d+n}) := \bigg\{ &\varphi\in C_b^{1,2}(\Lambda^{t,d+n})\colon\text{there exists $\textup{\text{\footnotesize H}}\in{\mathcal{H}}^{t,d+n}$ such that} \\
&\quad 0 = (\varphi - \hat u^{t,\omega^{d+n}})(t,0) = \underline{\mathcal{S}}_t^L \Big[\big(\varphi - \hat u^{t,\omega^{d+n}}\big)\big(\cdot\wedge\textup{\text{\footnotesize H}},B^{t,d+n}\big)\Big] \bigg\}. \notag
\end{align}
Let $\hat u \in \overline{\mathcal{U}}^{d+n}$. For every $L>0$ and $(t,\omega^{d+n})\in[0,T)\times\Omega^{d+n}$, define:
\begin{align*}
\overline{\mathcal{A}}^L \hat u(t,\omega^{d+n}) := \bigg\{ &\varphi\in C_b^{1,2}(\Lambda^{t,d+n})\colon\text{there exists $\textup{\text{\footnotesize H}}\in{\mathcal{H}}^{t,d+n}$ such that} \\
&\quad 0 = (\varphi - \hat u^{t,\omega^{d+n}})(t,0) = \overline{\mathcal{S}}_t^L \Big[\big(\varphi - \hat u^{t,\omega^{d+n}}\big)\big(\cdot\wedge\textup{\text{\footnotesize H}},B^{t,d+n}\big)\Big] \bigg\}.
\end{align*}
\end{Definition}

\begin{Definition}
\label{ViscSol}
\
\begin{itemize}
\item[(i)] \textup{Viscosity $L$-subsolution (resp., $L$-supersolution):} Let $\hat u \in \underline{\mathcal{U}}^{d+n}$ $($resp., $\hat u \in \overline{\mathcal{U}}^{d+n}$$)$. For any $L>0$, we say $\hat u$ is a viscosity $L$-subsolution $($resp., $L$-supersolution$)$ to PPDE \reff{PPDE} if
\begin{align*}
&-\partial_t \varphi(t,0) - (\hat{\mathcal{L}}^{t,\omega^{d+n}} \varphi)(t,0) \\
&- f\big(t,\omega^d,\omega_t^n,\hat u(t,\omega^{d+n}),\partial_{\omega^d} \varphi(t,0) + \partial_{\omega^n} \varphi(t,0) \sigma(t,\omega^d,\omega_t^n)\big)\leq 0, \qquad \text{$($resp., $\geq0$$)$},
\end{align*}
for any $(t,\omega^{d+n})\in[0,T)\times\Omega^{d+n}$ and any $\varphi\in\underline{\mathcal{A}}^L \hat u(t,\omega^{d+n})$ $($resp., $\varphi\in\overline{\mathcal{A}}^L \hat u(t,\omega^{d+n})$$)$, where
\begin{align*}
(\hat{\mathcal{L}}^{t,\omega^{d+n}} \varphi)(t,0)&:= \big\langle b^{t,\omega^d}(t,0,\omega_t^n),\partial_{\omega^n} \varphi(t,0)\big\rangle + \frac{1}{2}\textup{tr}\Big(\partial_{\omega^d\omega^d}^2 \varphi(t,0)(\sigma\sigma^T)^{t,\omega^d}(t,0,\omega_t^n) \\
&\quad + 2\partial_{\omega^d\omega^n}^2 \varphi(t,0)\sigma^{t,\omega^d}(t,0,\omega_t^n) + \partial_{\omega^d\omega^d}^2 \varphi(t,0)\Big).
\end{align*}
\item[(ii)] \textup{Viscosity subsolution (resp., supersolution):} Let $\hat u \in \underline{\mathcal{U}}^{d+n}$ $($resp., $\hat u \in \overline{\mathcal{U}}^{d+n}$$)$. We say $\hat u$ is a viscosity subsolution $($resp., supersolution$)$ to PPDE \reff{PPDE} if $\hat u$ is a viscosity $L$-subsolution $($resp., supersolution$)$ to PPDE \reff{PPDE} for some $L>0$.
\item[(iii)] \textup{Viscosity solution:} Let $\hat u \in UC_b(\Lambda^{d+n})$. We say $\hat u$ is a viscosity solution to PPDE \reff{PPDE} if it is both a viscosity subsolution and a viscosity supersolution.
\end{itemize}
\end{Definition}

\begin{Remark}
\label{R:viscdef}
{\rm We refer the reader to \cite{etzI} (see also \cite{ektz}) for several remarks about the above definition of viscosity solution of the path-dependent PDE \reff{PPDE}. Here we focus only on a point which will be useful later on, namely that the viscosity property is a \emph{local property}. More precisely, for any $(t,\omega^{d+n})\in[0,T)\times\Omega^{d+n}$ and any $\varepsilon>0$, define
\begin{equation}
\label{Heps}
\textup{\text{\footnotesize H}}_\varepsilon^t := \inf\big\{s>t\colon|B_s^{t,d+n}|\geq\varepsilon\big\}\wedge(t+\varepsilon) \qquad \text{and} \qquad \textup{\text{\footnotesize H}}_\varepsilon := \textup{\text{\footnotesize H}}_\varepsilon^0.
\end{equation}
Note that $\textup{\text{\footnotesize H}}_\varepsilon^t\in{\mathcal{H}}^{t,d+n}$. Then, to check the viscosity property of $\hat u$ at $(t,\omega^{d+n})$, it suffices to know the value of $\hat u^{t,\omega^{d+n}}$ on $[t,\textup{\text{\footnotesize H}}_\varepsilon^t]$ for an arbitrary small $\varepsilon>0$. In particular, for any $\varphi\in\underline{\mathcal{A}}^L \hat u(t,\omega^{d+n})$ with corresponding $\textup{\text{\footnotesize H}}\in{\mathcal{H}}^{t,d+n}$, we have $\textup{\text{\footnotesize H}}_\varepsilon^t \leq \textup{\text{\footnotesize H}}$ when $\varepsilon$ is small enough.}
\end{Remark}

\begin{Theorem}
\label{ThmEx}
Let assumptions \textup{\textbf{(HFC)}} and \textup{\textbf{(HBC)}} hold. Then the map $\hat v$ defined in \reff{hatv} belongs to $UC_b(\Lambda^{d+n})$ and is a viscosity solution to PPDE \reff{PPDE}.
\end{Theorem}
\emph{Proof.}
Thanks to Proposition~\ref{UC} and identification \reff{vhatv} we have $\hat v\in UC_b(\Lambda^{d+n})$. Now we prove the viscosity property.

We only show that $\hat v$ is a viscosity subsolution, the supersolution property is proved analogously. Without loss of generality, we prove the viscosity subsolution property at $(t,\omega^{d+n})=(0,0)$ and, for notational simplicity, we omit the superscript $^{0,0}$ in the rest of the proof. Moreover, in what follows we shall denote by $\underline L>0$ a generic positive constant depending only on the Lipschitz constant of $f$ with respect to $(y,z)$, which may vary from line to line.

Let us proceed by contradiction. Then, for any $L\geq\underline L$, there exists $\varphi\in\underline{\mathcal{A}}^L \hat v(0,0)$ such that
\[
-\big\{\partial_t \varphi + \hat{\mathcal{L}}\varphi + f(\cdot,\hat v,\partial_{\omega^d} \varphi + \partial_{\omega^n} \varphi\,\sigma)\big\}(0,0)=:c>0.
\]
Let $\text{\footnotesize H}\in{\mathcal{H}}^{d+n}$ be the hitting time corresponding to $\varphi$ in \reff{AL}. Let us denote simply by $(\hat X,\hat Y,\hat Z) := (\hat X^{0,0},\hat Y^{0,0},\hat Z^{0,0})$ the solution to the FBSDE \reff{FBSDEd+n} starting at time $t=0$ from $\omega^{d+n}=0$. Let $\varepsilon>0$ and define
\begin{equation}
\label{Heps'}
\textup{\text{\footnotesize H}}_\varepsilon' := \inf\big\{t>0\colon|(B_s^d,\hat X_s)|\geq\varepsilon\big\} \wedge \textup{\text{\footnotesize H}} \wedge \textup{\text{\footnotesize H}}_\varepsilon
\end{equation}
where $\textup{\text{\footnotesize H}}_\varepsilon$ is given by \reff{Heps}. Since $\varphi \in C_b^{1,2}(\Lambda^{d+n})$, by the uniform continuity of $\hat v$, $b$, $\sigma$, and $f$, we may assume $\varepsilon$ small enough such that
\begin{equation}
\label{c/2}
-\big\{\partial_t \varphi + \hat{\mathcal{L}}\varphi + f(\cdot,\hat v,\partial_{\omega^d} \varphi + \partial_{\omega^n} \varphi\,\sigma)\big\}(t,\omega^{d+n}) \geq \frac{c}{2}>0 \qquad\quad t\in[0,\text{\footnotesize H}_\varepsilon'].
\end{equation}
Now we observe that the following dynamic programming principle holds
\begin{equation}
\label{DPP}
\hat Y_0 = \hat v(\tau,B^d,\hat X) + \int_0^\tau f(s,B^d,\hat X_s,\hat Y_s,\hat Z_s)ds - \int_0^\tau \hat Z_s d B_s^d
\end{equation}
for any $\tau\in{\mathcal{T}}^{d+n}$, $\P_0^{d+n}$-a.s.. As a matter of fact, for any $t\in[0,T]$, by taking the conditional expectation $\E_{t}^{\P_0^{d+n}}$ and using the properties of regular conditional probability distributions, we deduce that $\hat Y$ is the solution to the following BSDE on $[t,T]$: for $\P_0^{d+n}$-a.e. $B^{d+n}$,
\[
\hat Y_s = g^{t,B^d}(B^{t,d},\hat X_T^{t,(B^d,\hat X)}) + \int_s^T f^{t,B^{d}}(r,B^{t,d},\hat X_r^{t,(B^d,\hat X)},\hat Y_r,\hat Z_r)dr - \int_s^T \hat Z_r dB_r^{t,d}
\]
for all $t \leq s \leq T$, $\P_0^{t,d+n}$-a.s.. Then $\hat Y_t$ $=$ $\hat Y_t^{t,(B^d,\hat X)}$, $\P_0^{d+n}$-a.s. and therefore $\hat v(t,B^d,\hat X)$ $=$ $\hat Y_t$, $\P_0^{d+n}$-a.s.. Thanks to the regularity in $t$ of $\hat v$, we obtain \reff{DPP}.

Then, from the functional It\^o's formula and Remark \ref{BetaSymm}, we have
{\allowdisplaybreaks
\begin{align*}
0 &= (\varphi - \hat v)(\text{\footnotesize H}_\varepsilon',B^d,\hat X) - \int_0^{\text{\footnotesize H}_\varepsilon'} \big(\partial_{\omega^d} \varphi(s,B^d,\hat X) + \partial_{\omega^n} \varphi(s,B^d,\hat X) \sigma(s,B^d,\hat X_s) - \hat Z_s\big) dB_s^d \\
&\quad- \int_0^{\text{\footnotesize H}_\varepsilon'} \big(\partial_t \varphi(s,B^d,\hat X) + \hat{\mathcal{L}}\varphi(s,B^d,\hat X) + f(s,B^d,\hat X_s,\hat v(s,B^d,\hat X),\hat Z_s)\big)ds.
\end{align*}}
Now, from \reff{Heps'} and \reff{c/2}, we get
{\allowdisplaybreaks
\begin{align*}
0 &\geq (\varphi - \hat v)(\text{\footnotesize H}_\varepsilon',B^d,\hat X) - \int_0^{\text{\footnotesize H}_\varepsilon'} \big(\partial_{\omega^d} \varphi(s,B^d,\hat X) + \partial_{\omega^n} \varphi(s,B^d,\hat X) \sigma(s,B^d,\hat X_s) - \hat Z_s\big) dB_s^d \\
&\quad+ \int_0^{\text{\footnotesize H}_\varepsilon'} \Big(\frac{c}{2} - f(s,B^d,\hat X_s,\hat v(s,B^d,\hat X),\partial_{\omega^d} \varphi(s,B^d,\hat X) + \partial_{\omega^n} \varphi(s,B^d,\hat X)\sigma(s,B^d,\hat X_s)) \\
&\quad+ f(s,B^d,\hat X_s,\hat v(s,B^d,\hat X),\hat Z_s)\Big)ds.
\end{align*}}
Therefore, there exists an $\R^d$-valued $\F^{d+n}$-progressively measurable stochastic process $\alpha$, with $|\alpha|\leq\underline L$, such that
{\allowdisplaybreaks
\begin{align*}
0 &\geq (\varphi - \hat v)(\text{\footnotesize H}_\varepsilon',B^d,\hat X) - \int_0^{\text{\footnotesize H}_\varepsilon'} \big(\partial_{\omega^d} \varphi(s,B^d,\hat X) + \partial_{\omega^n} \varphi(s,B^d,\hat X)\sigma(s,B^d,\hat X_s) - \hat Z_s\big) dB_s^d \\
&\quad+ \frac{c}{2}\text{\footnotesize H}_\varepsilon' + \int_0^{\text{\footnotesize H}_\varepsilon'} \big(\partial_{\omega^d} \varphi(s,B^d,\hat X) + \partial_{\omega^n} \varphi(s,B^d,\hat X)\sigma(s,B^d,\hat X_s) - \hat Z_s\big)\alpha_s ds \\
&= (\varphi - \hat v)(\text{\footnotesize H}_\varepsilon',B^d,\hat X) + \frac{c}{2}\text{\footnotesize H}_\varepsilon' \\
&\quad - \int_0^{\text{\footnotesize H}_\varepsilon'} \big(\partial_{\omega^d} \varphi(s,B^d,\hat X) + \partial_{\omega^n} \varphi(s,B^d,\hat X)\sigma(s,B^d,\hat X_s) - \hat Z_s\big) (dB_s^d - \alpha_s ds).
\end{align*}}
Now, consider the probability measure $\P^{d+n}$ equivalent to $\P_0^{d+n}$ on $(\Omega^{d+n},\mathcal{F}_T^{d+n})$ with Radon-Nikodym density:
\[
\frac{d\P^{d+n}}{d\P_0^{d+n}}\bigg|_{\mathcal{F}_t^{d+n}} = M_t := \exp\bigg(\int_0^t \alpha_s dB_s^d - \frac{1}{2}\int_0^t\alpha_s ds\bigg), \qquad 0 \leq t \leq T,\;\P_0^{d+n}\text{-a.s.}
\]
Then
\[
\E^{\P^{d+n}}\big[(\varphi - \hat v)(\text{\footnotesize H}_\varepsilon',B^d,\hat X)\big] < 0.
\]
This is a contradiction, since $\varphi \in \underline{\mathcal{A}}^L\hat v(0,0)$.
\ep

\vspace{2mm}

Finally, we provide a uniqueness result for the path-dependent PDE \reff{PPDE}. Up to now, the finest comparison theorem which may be used in our framework is Theorem 3.4 in \cite{etzII}. Note however that it requires strong assumptions on the forward coefficients, as stated in the following theorem.

\begin{Theorem}
\label{CompThm}
Let assumptions \textup{\textbf{(HFC)}} and \textup{\textbf{(HBC)}} hold. Suppose also that $b$ and $\sigma$ are constants and $\sigma\sigma^T$ is positive definite.
\begin{enumerate}
\item Consider a bounded viscosity subsolution $\hat u^1$ and a bounded viscosity supersolution $\hat u^2$ to PPDE \reff{PPDE}. If $\hat u^1(T,\omega^{d+n}) \leq g(\omega^d,\omega_t^n) \leq \hat u^2(T,\omega^{d+n})$, for all $\omega^{d+n}\in\Omega^{d+n}$, then $\hat u^1 \leq \hat u^2$ on $\Lambda^{d+n}$.
\item The map $\hat v$ given by \reff{hatv} is the unique viscosity solution in $UC_b(\Lambda^{d+n})$ to the PPDE \reff{PPDE} which satisfies the terminal condition \reff{termcondPPDE}.
\end{enumerate}
\end{Theorem}
\textbf{Proof.}
Theorem \reff{CompThm} follows directly from Theorem 3.4 in \cite{etzII}. Indeed, it is easy to show that all the hypotheses required by Theorem 3.4 in \cite{etzII} are satisfied. In particular, assumption 2.8 in \cite{etzII} is satisfied, since $b$ and $\sigma$ are constants. Furthermore, assumption 3.1 in \cite{etzII} is a consequence of Proposition 7.2 in \cite{etzII} and the fact that $\sigma\sigma^T$ is positive definite.
\ep

\begin{Remark}
\label{R:Enlarge}
Enlarging the space. {\rm We would like to stress that the method of enlarging the space presented in this section may also be used in different contexts. As an example, for any $(t,\omega)\in\Lambda^n$, let us consider the following FBSDE:
\begin{equation}
\label{FBSDE2}
\begin{cases}
X_s = \omega_t + \int_t^s \tilde b^{t,\omega}(r,X)dr + \int_t^s \tilde \sigma^{t,\omega}(r,X)dB_r^{t,d} \\
Y_s = \tilde g^{t,\omega}(X) + \int_s^T \tilde f^{t,\omega}(r,X,Y_r,Z_r)dr - \int_s^T Z_rdB_r^{t,d}
\end{cases}
\end{equation}
for all $t \leq s \leq T$, $\P_0^{t,d}$-a.s., where $\tilde b\colon\Lambda^n\rightarrow\R$, $\tilde \sigma\colon\Lambda^n\rightarrow\R^{n\times d}$, $\tilde g\colon\Omega^n\rightarrow\R$, and $\tilde f\colon\Lambda^n\times\R\times\R^d\rightarrow\R$ satisfy assumptions \textbf{(HFC)} and \textbf{(HBC)}, adapted to the present context in an obvious way. Define
\begin{equation}
\label{v2}
\tilde v(t,\omega) := Y_t^{t,\omega}, \qquad (t,\omega)\in\Lambda^n.
\end{equation}
The FBSDE \reff{FBSDE2} has also been considered in \cite{etzI,h-ltt13} for the case $n=d$ and $\tilde\sigma\tilde\sigma^T$ strictly positive. In contrast, without imposing these conditions, we may proceed as at the beginning of this section and rewrite \reff{FBSDE2} on $(\Omega^{d+n},\mathcal{F}_T^{d+n},\P_0^{d+n})$. We may then define
\begin{equation}
\label{hatv2}
\hat{\tilde{v}}(t,\omega^{d+n}) :=  \tilde v(t,\omega^n), \qquad (t,\omega^{d+n})\in\Lambda^{d+n}.
\end{equation}
The path-dependent PDE associated to \reff{FBSDE2} is given by:
\begin{align}
\label{PPDE2}
-\partial_t \hat{\tilde{u}}(t,\omega^{d+n}) &- (\hat{\tilde{\mathcal{L}}} \hat{\tilde u})(t,\omega^{d+n}) \\
&- \tilde f\big(t,\omega^n,\hat{\tilde u}(t,\omega^{d+n}),\partial_{\omega^d} \hat{\tilde u}(t,\omega^{d+n}) + \partial_{\omega^n} \hat{\tilde u}(t,\omega^{d+n})\tilde\sigma(t,\omega^n)\big) = 0, \notag
\end{align}
for all $(t,\omega^{d+n})\in[0,T)\times\Omega^{d+n}$, with the terminal condition:
\begin{equation}
\label{termcondPPDE2}
\hat{\tilde u}(T,\omega^{d+n}) = \tilde g(\omega^n), \qquad\qquad \omega^{d+n}\in\Omega^{d+n},
\end{equation}
where
\begin{align*}
(\hat{\tilde {\mathcal{L}}} \hat{\tilde u})(t,\omega^{d+n})&:= \big\langle \tilde b(t,\omega^n),\partial_{\omega^n} \hat{\tilde u}(t,\omega^{d+n})\big\rangle + \frac{1}{2}\text{tr}\Big(\partial_{\omega^n\omega^n}^2 \hat{\tilde u}(t,\omega^{d+n}) \tilde\sigma(t,\omega^n) \tilde\sigma^T(t,\omega^n) \\
&\quad + 2\partial_{\omega^d\omega^n}^2 \hat{\tilde u}(t,\omega^{d+n}) \tilde\sigma(t,\omega^n) + \partial_{\omega^d\omega^d}^2 \hat{\tilde u}(t,\omega^{d+n})\Big),
\end{align*}
for all $(t,\omega^{d+n})\in[0,T)\times\Omega^{d+n}$.

Proceeding as in Theorem \ref{ThmEx}, we may prove that \reff{hatv2} is a viscosity solution to \reff{PPDE2}-\reff{termcondPPDE2}. However, since $\hat{\tilde{v}}=\hat{\tilde{v}}(t,\omega^d,\omega^n)$ is constant with respect to $\omega^d$, then it is simple to show that $\hat{\tilde{v}}$ (and hence $\tilde{v}$) is a viscosity solution to the following path-dependent PDE:
\begin{equation}
\label{PPDE3}
-\partial_t \tilde u(t,\omega) - (\tilde{\mathcal{L}} \tilde u)(t,\omega) - \tilde f\big(t,\omega,\tilde u(t,\omega),\partial_{\omega} \tilde u(t,\omega)\tilde\sigma(t,\omega)\big) = 0, \notag
\end{equation}
for all $(t,\omega)\in[0,T)\times\Omega^n$, with the terminal condition:
\begin{equation}
\label{termcondPPDE3}
\tilde u(T,\omega) = \tilde g(\omega), \qquad\qquad \omega\in\Omega^n,
\end{equation}
where
\[
(\tilde {\mathcal{L}} \tilde u)(t,\omega):= \big\langle \tilde b(t,\omega),\partial_{\omega} \tilde u(t,\omega)\big\rangle + \frac{1}{2}\text{tr}\Big(\partial_{\omega\omega}^2 \tilde u(t,\omega) \tilde\sigma(t,\omega) \tilde\sigma^T(t,\omega)\Big),
\]
for all $(t,\omega)\in[0,T)\times\Omega^n$.}
\end{Remark}

\textbf{Acknowledgements.}
The author would like to take this opportunity to thank Professor Marco Fuhrman and Professor Nizar Touzi for helpful discussions and suggestions related to this work.

\end{document}